\documentclass[graybox]{svmult}

\usepackage{type1cm}        
\usepackage{makeidx}         
\usepackage{graphicx}        
\usepackage{multicol}        
\usepackage[bottom]{footmisc}

\usepackage{newtxtext}       %
\usepackage[varvw]{newtxmath}       

\makeindex             

\usepackage[pdftex]{hyperref}
\theoremstyle{definition}
\newtheorem{ex}{Example}
\DeclareMathOperator{\sech}{sech}
\DeclareMathOperator{\li}{Li}

\newcommand{\rmd}{\mathrm}
\newcommand{\disp}{\displaystyle}
\newcommand{\zp}{\zeta^{\prime}}

\newcommand{\G}{\textbf{\textup{G}}}

\begin{document}

\title*{A new approach to evaluating Malmsten's integral and related integrals}
\author{Abdulhafeez A. Abdulsalam}
\institute{Abdulhafeez A. Abdulsalam \at University of Ibadan, Department of Mathematics, University of Ibadan, Ibadan, Oyo, Nigeria, \email{hafeez147258369@gmail.com}}
%
%
\maketitle

\abstract{This paper discusses generalizations of logarithmic and hyperbolic integrals. A new proof for an integral presented by Vardi and several other integrals in relation to known mathematical constants are discovered. We introduce the signed generalized Stirling polynomials of the first kind from the generalized Stirling polynomials of the first kind, and we give new expressions for the signed generalized Stirling polynomials of the first kind in terms of the Stirling cycle numbers and complete Bell polynomials. We establish the role of the signed generalized Stirling polynomials of the first kind and complete Bell polynomials in generalizing Malmsten's integral for all natural powers of the hyperbolic secant function, and we derive a reduction formula for the integral sequence. We give expressions for new integral sequences, which possess similar properties with Malmsten's integral, in terms of the signed generalized Stirling polynomials of the first kind, and we discover identities and a functional equation for the signed generalized Stirling polynomials of the first kind.}

\section{Introduction}
Malmsten and colleagues \cite{bib13} treated a certain generalization of a logarithmic integral in 1842. Their exposition yields the integral
\begin{equation}
\int_0^\infty \ln{x}\sech{x}\,\rmd{d}x = \pi\ln\left(\frac{\Gamma\left(\frac{3}{4}\right)\sqrt{2\pi}}{\Gamma\left(\frac{1}{4}\right)}\right). \label{he2r}
\end{equation}
In 1972, Glasser \cite[(A14)]{bib28} evaluated the derivative of the Dirichlet beta function at 0. We discover that this value, combined with the derivative of the Dirichlet beta function at 1,  immediately yields \eqref{he2r}.  In a remarkable paper which appeared in the American Mathematical Monthly in 1988, Vardi \cite{bib17} treated several interesting logarithmic integrals found in Gradshteyn and Ryzhik \cite[(4.229.4)]{bib25}. Vardi proposed a method of proof based on the use of the Dirichlet $L$-function and Blagouchine showed that his result in the hyperbolic form of  \eqref{he2r} follows immediately from Malmsten's result \cite[p.~32]{bib1}. Blagouchine did not mention Glasser's result in his impressive exposition. Perhaps, it was not known to Blagouchine that Glasser had evaluated the derivative of the Dirichlet beta function at 0 in the early 1970's; almost two decades before Vardi evaluated \eqref{he2r}. In Moll's treatise on the proofs of special integrals of Gradshteyn and Ryzhik, published months after Blagouchine's paper, Vardi's proof for \eqref{he2r} was presented but neither Blagouchine's proof \cite[\S3.4.1]{bib1}, nor the proof to be presented in this exposition were discussed \cite[p.~2]{bib7}. In this paper, we also evaluate \eqref{he2r} by a method simpler than the methods explored by the aforementioned authors.\\

We extend \eqref{he2r} to all natural powers of $\sech{x}$ as follows:
\begin{equation}
\Delta_n := \int_0^\infty \ln{x}\sech^n{x} \, \rmd{d}x, \label{he1r}
\end{equation}
where $n \in \mathbb{N}:=\{1,2,3,\ldots\}$. Blagouchine's method for evaluating \eqref{he1r} is entirely based on the knowledge of contour integration \cite[\S3.4.1]{bib1}, while the method for evaluating \eqref{he1r} in this work involves the use of the Lerch transcendent and the Hurwitz zeta function. To prove Vardi's result elementarily, we evaluate \eqref{he2r} with $\ln{x}$ replaced with $\ln{(x^2 + a^2)}$ for real values of $a$. The latter integral yields \eqref{he2r} at $a=0$. Subsequently, we use, for the first time, an integral found in V\u{a}lean \cite[\S 3.40]{bib9} to derive $\Delta_2$. We derive, by the use of two other related integral sequences, a reduction formula for \eqref{he1r} in terms of a double integral. And by the means of Ripon's generalization of the $\ln$-$\ln$ form of $\Delta_n$ \cite[(2.1)]{bib16}, we provide a generalization of \eqref{he1r}.\\

The outline of this article is as follows. In Section \ref{sec2}, we introduce the signed generalized Stirling polynomials of the first kind from the generalized Stirling polynomials of the first kind, we give new expressions for the signed generalized Stirling polynomials of the first kind in terms of the Stirling cycle numbers and complete Bell polynomials, and in a lemma, we prove, for the first time, a useful identity for the Stirling polynomials. This identity allows us to simplify expressions for more definite integral sequences which possess similar properties with $\Delta_n$. In Section \ref{secvardi}, we demonstrate that Glasser's derivative of the Dirichlet beta function at 0, together with the derivative of the Dirichlet beta function at 1, yields \eqref{he2r}. We also present a new elementary proof for integral \eqref{he2r} in the hyperbolic form. The proof of \eqref{he2r} presented in this work appears so far to be the simplest in the literature. In Section \ref{section4}, a new reduction formula for $\Delta_n$ is derived through integration by parts. We derive a new generalization of $\Delta_n$ in terms of the signed generalized Stirling polynomials of the first kind. This generalization is used to derive expressions for more definite integral sequences whose integrands are combinations of hyperbolic functions. We discover other identities and a functional equation for the signed generalized Stirling polynomials of the first kind. We reevaluate an integral presented in V\u{a}lean's book \cite[\S 3.40]{bib9} by a new approach. This reevaluation leads to the rediscovery of the value of $\Delta_2$. It also leads to a discovery of new integrals in Propositions \ref{prop4} and \ref{prop4w}, which discusses the relationship between known mathematical constants and integrals of certain combinations of transcendental functions. The two Computer Algebra System (CAS) softwares used to verify our results are the \textsf{Maple~2022} and \textsf{Mathematica~13}.\\

In summary, this paper establishes the following results as follows. Theorem \ref{newthm} and Lemmas \ref{lem}, \ref{lem2} have been used to write the values of $\Delta_n$ from $n=2$ to $6$. Blagouchine provided a generalization of $\Delta_1$ by contour integration, but we have given a solution using an elementary method. Remarks \ref{rmk2.2}, \ref{rmk4.1}, Propositions \ref{ssect3.7}, \ref{prop2}, \ref{prop3}, \ref{prop4}, \ref{prop4w}, Lemmas \ref{frlmk}, \ref{lem0}, \ref{lemconj}, and Theorems \ref{newthm}, \ref{conjecture} do not seem to be presented elsewhere. Lemma \ref{lem} was derived from previously existing results and Example \ref{scex} can be found in Blagouchine's paper. Lemma \ref{lem2} is not new since a generalization of it can be found in Blagouchine's paper.  We only point out that Lemma \ref{lem2} may be derived from Theorem \ref{newthm} and \eqref{he2r}, because it is necessary to simplify our expressions for the integral sequences. Since a generalization of Proposition \ref{propf} can be found in Blagouchine's paper, Proposition \ref{propf} and Remark \ref{rmk:lamb} are not new. We offer a new proof of Proposition \ref{propf} and derive  Remark \ref{rmk:lamb} from  \eqref{32}. We have not come across any material discussing $\lambda_n$ and $\delta_n$. However, we do point out that \textsf{Mathematica~13} gives closed form expressions for $\lambda_n$ for each $n$, and does not give any expression for $\delta_n$. The integral sequences $\chi_n$ and $\delta_n$ do not seem to have appeared elsewhere in the literature.\medskip
\section{Notations and Definitions}\label{sec2}
Throughout this manuscript, the following abbreviated notations are used: $\mathbb{N}_0 = \mathbb{N} \cup \{0\}$, $\gamma = 0.5772156649...$ represents Euler's constant, $\mathrm{e}=2.71828182845...$ represents Euler's number, $\G = 0.9159655941...$ represents Catalan's constant, $H_n$ represents the $n$-th harmonic number, and $B_n(x)$ represents the $n$-th Bernoulli polynomial \cite[\S 24.2]{bib23}. We denote, respectively, the gamma, digamma, and polygamma functions of argument $z$ with $\Gamma(z)$, $\psi(z)$ and $\psi_n(z)$, where $n \in \mathbb{N}$ \cite[(5.2.1), (5.2.2), \S 5.15]{bib23}. The Lerch transcendent is defined as \cite[(24.14.1)]{bib23}:
$$\Phi(z, s, a) = \sum_{n=0}^{\infty} \frac{z^n}{(n + a)^s}, \quad \lvert z\rvert \leq 1, \Re\,s > 1, a \not\in -\mathbb{N}_0.$$
The polylogarithm function \cite[\S 25.12(ii)]{bib23} $\li_s(z)$, is defined as: $\li_s(z) = z\Phi(z, s, 1)$. The Riemann zeta function \cite[\S 25.2]{bib23} and the Hurwitz zeta function \cite[\S 25.11]{bib23} are, respectively, defined as:
$$\zeta(s) = \sum_{n=1}^{\infty} \frac{1}{n^s}, \quad \zeta(s, z) = \sum_{n=0}^{\infty} \frac{1}{(n+z)^s},$$
where $z \not\in -\mathbb{N}_0$, $\Re\,s > 1$. The domain $\Re\,s > 1$ can be extended to $s \in \mathbb{C}\setminus \{1\}$ through analytic continuation, using for instance, the Hermite integral representation for the Hurwitz zeta function \cite[(25.11.29)]{bib23}. We denote the $n$-th derivative of $\zeta(s,z)$  with respect to $s$ as $\zeta^{(n)}(s,z)$, namely $\zeta^{(n)}(s,z) := \frac{\rmd{d}^{n}}{\rmd{d}s^{n}}\zeta(s,z).$ The Dirichlet beta function is defined in \cite[(3)]{bib28} as:
$$\beta(s) =  \sum_{k=0}^{\infty} \frac{(-1)^k}{(2k+1)^s}, \quad \Re\,s >0.$$
An integral representation for $\beta(s)$, which provides analytic continuation of the domain of $\beta(s)$ to $\Re\,s<0$, can be found in Glasser's paper \cite[(A13)]{bib28}. The signed Stirling numbers of the first kind are denoted by $s(n, k)$, and the Stirling cycle numbers or unsigned Stirling numbers of the first kind are denoted by ${n \brack k}$. The signed Stirling numbers of the first kind are related to the Stirling cycle numbers by ${n \brack k} = (-1)^{k-n}s(n, k) = |s(n, k)|$ \cite[\S26.1]{bib23}. In combinatorics, the signed Stirling number of the first kind is $(-1)^{n-k}$ times the Stirling cycle number, which is the number of permutations of $\{1,2,3,\ldots,n\}$ with exactly $k$ disjoint cycles. The Stirling cycle numbers are defined by the generating functions \cite[(20.1), (20.2)]{bib27}
\begin{equation}\label{modff}
(x-n+1)_n= \sum_{k=0}^n (-1)^{n-k} {n \brack k}  x^k, \quad \sum_{n=k}^{\infty} (-1)^{n-k} {n \brack k} \frac{z^n}{n!} = \frac{1}{k!}\ln^k(1+z),
\end{equation}
where $(a)_n$ is the rising factorial defined by $(a)_n := a(a+1)(a+2)\cdots(a+ n-1)$. It is worth noting that we slightly modified \eqref{modff} from their sources by using the property ${n \brack k} = (-1)^{k-n}s(n, k)$. Ripon \cite[(4.34)]{bib16} defined the polynomials $P_k\left(m, x\right)$ by the generating function
\begin{equation}\label{gen2}
(j+1)_m = \sum_{k=0}^m (-1)^k P_k\left(m, x\right) \left(j + x\right)^{m-k},
\end{equation}
We shall refer to the polynomials $P_k\left(m, x\right)$ as the signed generalized Stirling polynomials of the first kind. The complete Bell polynomials are defined by the generating function \cite[(8)]{bib31}
\begin{equation}
\exp\left(\sum_{r=1}^{\infty} x_r \frac{t^r}{r!}\right) = \sum_{n=0}^{\infty} B_{n}(x_1,x_2,\ldots,x_n) \frac{t^n}{n!},
\end{equation}
and are given by \cite[(9)]{bib31}
\begin{equation}
B_{n}(x_1,x_2,\ldots,x_n) = \sum_{\pi(n)} \frac{n!}{\prod_{r=1}^n k_r!}\prod_{r=1}^{n} \left(\frac{x_r}{r!}\right)^{k_r},
\end{equation}
where the sum is taken over all partitions $\pi(n)$ of $n$, that is, over all $k_1, k_2, k_3, \ldots, k_n \in \mathbb{N}_0$ satisfying $\sum_{r=1}^n rk_r = n$. The generalized harmonic number, $H_{n}^{(m)}$, is defined as:
$$H_{n}^{(m)} := \sum_{k=1}^n \frac{1}{k^m} = \frac{(-1)^m}{(m-1)!}\left(\psi_{m-1}(1) - \psi_{m-1}(n+1)\right),\quad m \in \mathbb{N},\, n\in \mathbb{N}_0.$$
We now define the definite integral sequences $\lambda_n$, $\delta_n$, $\chi_n$, and provide relationships between $\Delta_n$, $\lambda_n$, $\delta_n$ and $\chi_n$.
\begin{definition}\label{def1} Let $n \in \mathbb{N}$. We define $\lambda_n$, $\delta_n$, and $\chi_n$ as follows
\begin{align*}
&\lambda_n := \int_0^{\infty}\frac{\tanh{x}\sech^n{x}}{x}\,\rmd{d}x, \\\\
&\delta_n := \int_0^{\infty} \frac{1 - \sech{x}}{x^2}\sech^n{x}\, \rmd{d}x, \\\\
&\chi_{n} := \int_0^{\infty}\frac{\sech{x} - \sech^n{x}}{x^2}\,\rmd{d}x.
\end{align*}
\end{definition}
\begin{remark}
The relationship between $\Delta_n$ and $\lambda_n$ is:
\begin{equation}\label{eqn1.5}
\Delta_n = \frac{n-2}{n-1}\Delta_{n-2} - \frac{1}{n-1}\lambda_{n-2}, \quad n \geq 3,
\end{equation}
and the relationship between $\lambda_n$ and $\delta_{n}$ is:
\begin{equation}\label{eqn1.6}
\lambda_n = \frac{1}{n} \delta_{n-1} + \frac{n-1}{n}\lambda_{n-1}, \quad n \geq 2.
\end{equation}
These follow directly from integration by parts. After some iterations, we obtain from \eqref{eqn1.6}
\begin{equation}\label{red2}
\lambda_n = \frac{1}{n}\sum_{k=1}^{n-1} \delta_{n-k} + \frac{1}{n}\lambda_1, \quad n \geq 2.
\end{equation}
We find that $\sum_{k=1}^{n-1} \delta_{n-k} = \sum_{k=1}^{n-1} \delta_{k} = \chi_n$, and $n\lambda_n = \chi_n + \lambda_1$. Using \eqref{eqn1.5} and \eqref{red2}, we derive the relationship between $\Delta_n$ and $\chi_n$ as:
\begin{equation}\label{red4}
\Delta_n = \frac{n-2}{n-1}\Delta_{n-2} - \frac{1}{(n-2)(n-1)}\chi_{n-2} - \frac{1}{(n-2)(n-1)}\lambda_1,\quad n \geq 3.
\end{equation} 
We find that $\lambda_1 = \delta_0$ and $\chi_2 = \delta_1$; the second being trivial.
\end{remark}
We now define three functions of $y \in \mathbb{R}$ that are useful in deriving the value of $\Delta_2$.
$$\kappa(y) := \int_0^\infty \tanh\left( \frac{\pi t}{2}\right)\left(\frac{1}{t} - \frac{t}{t^2 + y^2}\right) \mathrm{d}t, \quad \kappa_2(y)  := \int_0^\infty\left(\frac{\tanh\left(\frac{\pi t}{2}\right)}{t} - \frac{t}{t^2 + y^2}\right)\mathrm{d}t,$$
$$\kappa_1(y)  := \int_0^\infty\left(\frac{t}{t^2 + y^2} - \frac{t}{t^2 + y^2}\tanh\left(\frac{\pi t}{2}\right)\right) \mathrm{d}t,$$
where $y \neq 0$ for $\kappa_1(y)$ and $\kappa_2(y)$.
\begin{remark}
Ripon \cite[(4.34)]{bib16} gave the following expressions for the signed generalized Stirling polynomials of the first kind
$$P_k\left(m, x\right) := \binom{m}{k}x^k - \binom{m}{k}\frac{k(m+1)}{2}x^{k-1}  + \sum_{j=2}^k (-1)^j R(m, j) x^{k-j},$$
where $R(m, j)= 0$ for $m \leq j$, and
\begin{align*}
R(m, j) &= \frac{m!}{(m-j)!}B_{m-j}\left(H_m^{(1)}, -H_m^{(2)},\ldots, (-1)^{m-j-1}(m-j-1)!H_m^{(m-j)}\right)
\\&= \frac{1}{j!}B_{j}\left(H_m^{(-1)}, \ldots, (-1)^{j}(j-1)!H_m^{(-j)}\right),
\end{align*}
for $m > j$. A careful study of Ripon's expressions in relation to \eqref{gen2}, and an implementation in \textsf{Maple~2022} shows that Ripon's expressions fail to generate several values of the signed generalized Stirling polynomials. Upon this discovery, we are compelled to derive new working expressions for these polynomials. The expressions we derive in the sequel correctly generate values of the signed generalized Stirling polynomials. 
\end{remark}

\section{Useful lemmas}
In this section, we present some useful lemmas, of which the first two are new. In the forthcoming lemma, we derive, in terms of the Stirling cycle numbers and complete Bell polynomials, new expressions for the signed generalized Stirling polynomials of the first kind.
\begin{lemma}\label{frlmk}
Let $k, m \in \mathbb{N}_0$, $m  \geq k$ and $\Re\,x>0$. Then
\begin{align}
P_k(m, x) &= \sum_{j=0}^{k} (-1)^{j-k} x^{j} \binom{j+m-k}{m-k}{m+1 \brack j+m-k+1} \nonumber
\\&= \sum_{j=m-k}^{m} (-1)^{j-m} x^{j+k-m} \binom{j}{m-k} \frac{m!}{j!} \nonumber
\\&\qquad \times B_{j}\left(H_m^{(1)},-H_m^{(2)},\ldots, (-1)^{j-1} (j-1)! H_m^{(j)}\right).\label{lstr1}
\end{align}
\end{lemma}
\begin{proof}
The generalized Stirling polynomials of the first kind are defined by the generating function \cite[p.~7]{bib29} 
\begin{equation}\label{lstr3}
 \prod_{k=1}^{n-1} (k+x-z) = \sum_{k=0}^{n-1} 	P_{k, n}(z) x^k, \quad \Re\,z > 0.
\end{equation}
Adamchik \cite[(18)]{bib29} showed that
\begin{equation}\label{lstr2}
P_{k, n}(z) = \sum_{j=k+1}^n (-z)^{j-k-1}\binom{j-1}{k}{n \brack j}.
\end{equation}
By substituting $x = j+z$ and $n=m+1$ in \eqref{lstr3}, followed by a change of variable from $z$ to $x$, we derive
$$ \prod_{k=1}^{m} (k+j) = (j+1)_m  = \sum_{k=0}^{m} P_{k, m+1}(x) (j+x)^k.$$
We can express the sum on the right-hand side of \eqref{gen2} as
\begin{equation}\label{lstr4}
\sum_{k=0}^m (-1)^k P_k\left(m, x\right) \left(j + x\right)^{m-k} = \sum_{k=0}^m (-1)^{m-k} P_{m-k}(m, x) (j + x)^k.
\end{equation}
We discover from \eqref{lstr4} that the relationship between the polynomials $P_k(m, x)$ and the generalized Stirling polynomials of the first kind is: $P_{k, m+1}(x) = (-1)^{m-k} P_{m-k}(m, x)$. Hence \eqref{lstr1} follows from \eqref{lstr2}. We envisage the polynomials $P_k\left(m, x\right)$ as the signed generalized Stirling polynomials of the first kind, due to the relationship between $P_{k, n}(x)$ and $P_{k}(m, x)$. To prove the second representation, we note that ${n \brack r} = (-1)^{n-r} s(n, r)$, and with this, we easily deduce from \cite[(44.10)]{bib23}
\begin{equation}\label{hert}
{n+1 \brack r+1} = \frac{n!}{r!}B_r\left(H_n^{(1)}, -H_n^{(2)},\ldots,(-1)^{r-1}(r-1)!H_n^{(r)}\right).
\end{equation}
By making necessary substitutions in \eqref{hert}, we derive the second representation.
\end{proof}
In the next lemma, we introduce an important identity involving the Stirling polynomials, this identity was in fact, not mentioned in either Adamchik's or Ripon's paper. We also demonstrate the usefulness of the identity in evaluating $\Delta_n$.  
\begin{lemma}\label{lem0}Let $n, r\in \mathbb{N}$, $k =2r$ and $n \geq 2r$. Then
\begin{equation}\label{yufo}
P_{n-2r, n}\left(\frac{n}{2}\right) =  P_{2r-1}\left(n-1, \frac{n}{2}\right) = 0.
\end{equation}
\end{lemma}
\begin{proof}
The polynomials $P_{k, n}(x)$ have the alternative representation \cite[(19)]{bib29} 
$$P_{k, n}(x) = \lim_{y \to 0} \frac{(-1)^k}{k!}\frac{\rmd{d}^{n-1}}{\rmd{d}y^{n-1}}\frac{\ln^{k}(1-y)}{(1-y)^{1-x}}.$$
We derive by Leibniz theorem \cite{bib30} for even $n$, where $n \geq 2r$
\begin{align} 
&\frac{\rmd{d}^{n-1}}{\rmd{d}y^{n-1}} \ln^{n-2r}(1-y) \left(1 - y\right)^{\frac{n}{2}-1} \nonumber
\\&= \sum_{s=0}^{n-1}\sum_{k=1}^{n-s-1} \binom{n-1}{s}\binom{\frac{n}{2}-1}{s}\frac{(-1)^{k+s} s! (n-2r)! \phi(k, n-s-1)}{(n-2r- k)!\left(1 - y\right)^{\frac{n}{2}}}\ln^{n-2r-k}(1 -y),\label{idnd}
\end{align}
where $\phi: \mathbb{N}^2 \longrightarrow \mathbb{N}$. No expression is known for $\phi(k, n-s-1)$, but few values of $\phi(k, n-s-1)$ are: $\phi(1, 1) = 1$; $\phi(1,2) = 1,\,\phi(2,\,2) = 1$; $\phi(1,3) = 2,\,\phi(2,3) = 3,\,\phi(3,\,3) = 1$; $\phi(1,4) = 3!,\,\phi(2,4) = 11,\,\phi(3,4) = 6,\,\phi(4,4) = 1$; $\phi(1,5) = 4!,\,\phi(2,5) = 50,\,\phi(3,5)  = 35,\,\phi(4,5) = 10,\,\phi(4,5) = 1$. Taking limits on both sides of \eqref{idnd} as $y \to 0$ yields
$$\lim_{y\to0} \frac{\rmd{d}^{n-1}}{\rmd{d}y^{n-1}} \ln^{n-2r}(1-y) \left(1 - y\right)^{\frac{n}{2}-1} = 0.$$
We have for odd $n$, where $n > 2r$
\begin{align} 
&\lim_{y \to 0} \frac{\rmd{d}^{n-1}}{\rmd{d}y^{n-1}} \ln^{n-2r}(1-y) \left(1 - y\right)^{\frac{n}{2}-1} = \lim_{y,z \to 0}\frac{\rmd{d}^{n-1}}{\rmd{d}y^{n-1}} \frac{\rmd{d}^{n-2r}}{\rmd{d}z^{n-2r}} \left(1 - y\right)^{z+ \frac{n}{2}-1} \nonumber
\\&=  \lim_{z \to 0}\frac{\rmd{d}^{n-2r}}{\rmd{d} z^{n-2r}} \frac{(-1)^{n-1}\Gamma\left(z + \frac{n}{2}\right)}{\Gamma\left(z - \frac{n}{2}+1\right)} =  \lim_{z \to 0}\frac{\rmd{d}^{n-2r}}{\rmd{d} z^{n-2r}} (-1)^{n-1} \prod_{k=1}^{\frac{n-1}{2}} \left(z^2 - \frac{(2k -1)^2}{4}\right) \nonumber
\\& = \lim_{z \to 0} \sum_{k=0}^{\frac{n-1}{2}} \frac{(n-2k-1)! (-1)^{n-1}  R(k, n)}{(2r-2k-1)!}z^{2r-2k-1} = 0. \nonumber
\end{align}
Hence, \eqref{yufo} follows.
\end{proof}
\begin{remark}\label{rmk2.2}
Taking account of \eqref{gen2} and \eqref{yufo}, we find, for odd $n$
$$\prod_{k=1}^{\frac{n-1}{2}} \left(z^2 - \frac{(2k -1)^2}{4}\right) = \sum_{k=0}^{\frac{n-1}{2}} P_{2k}\left(n-1, \frac{n}{2}\right) z^{n-2k-1},$$
and so, $R(k, n) = P_{2k}\left(n-1, \frac{n}{2}\right)$. It can also be deduced from \eqref{lstr1} that $P_0(m, x) = P_0(0, x)  = 1$. 
\end{remark}
In the following lemma, we provide useful relations required in the simplification of the expressions for the integral sequences that will appear in later sections.
\begin{lemma}\label{lem}
Let $m,n \in \mathbb{N}$, $q \in \mathbb{N}_0$, $b > 0$, $\Re\,a > 1$, $x \in \mathbb{C}\setminus-\mathbb{N}_0$, $z \in (0, 1)$, and
$$\zeta^{(q)}\left(1, x\right) - \zeta^{(q)}\left(1, x+\frac{1}{2}\right) :=  \lim_{s\to1}\left[\zeta^{(q)}\left(s,  x\right) - \zeta^{(q)}\left(s,  x+\frac{1}{2}\right)\right].$$
Then
\begin{align}\label{fonl}
\zeta(a, b) - \zeta\left(a, b + \frac{1}{2}\right) = 2^a\Phi(-1, a, 2b),
\end{align}
\begin{align}\label{bernz}
\zeta(-n, x) - \zeta\left(-n, x + \frac{1}{2}\right) = \frac{1}{n+1}B_{n+1}\left(x + \frac{1}{2}\right) - \frac{1}{n+1}B_{n+1}(x),
\end{align}
\begin{align}\label{zpsi}
\zeta(n, x) - \zeta\left(n, x + \frac{1}{2}\right) = \frac{(-1)^{n}}{(n-1)!}\left(\psi_{n-1}\left(x\right) - \psi_{n-1}\left(x + \frac{1}{2}\right)\right),
\end{align}
\begin{align}\label{zzerox}
\zp(0, x) - \zp\left(0, x+ \frac{1}{2}\right) = \ln{\Gamma(x)} - \ln{\Gamma\left(x + \frac{1}{2}\right)},
\end{align}
\begin{equation}
\zp(-n, z) + (-1)^n \zp(-n, 1 - z) = \pi i \frac{B_{n+1}(z)}{n+1} + \frac{n!\mathrm{e}^{-\pi i n/2}}{(2\pi)^n}\li_{n+1}\left(\mathrm{e}^{2\pi iz}\right),\label{adamchik1}
\end{equation}
\begin{align}\label{prlf}
\zp\left(n, \frac{m}{2}\right) - \zp\left(n, \frac{m+1}{2}\right) =  &\frac{(-1)^{n}\ln{2}}{(n-1)!}\left(\psi_{n-1}\left(\frac{m}{2}\right) - \psi_{n-1}\left(\frac{m+1}{2}\right)\right) \nonumber
\\&+ q(m, n) + 2^n\left(-1\right)^{m} \left(1 - \delta_{m1}\right)\sum\limits_{k=1}^{m-1} \frac{(-1)^k\ln{k}}{k^n},
\end{align}
where 
\begin{equation}
q(m, n) := \begin{cases}2\ln{2}\left(-1\right)^{m-1}\zeta(n) + \left(2^{n} - 2\right)\zp(n),  &\textup{for}\,\,  n \neq 1,
\\\left(-1\right)^{m}\left( \ln^2{2} - 2\gamma\ln{2}\right), &\textup{for}\,\, n = 1,
\end{cases}
\end{equation}
and $\delta_{mn}$ is the Kronecker delta.
\begin{proof}
Separating the right side of \eqref{fonl} into odd and even parts
\begin{align*}
2^a\Phi(-1, a, 2b) &= \sum_{k=0}^{\infty} \frac{1}{\left(k + b\right)^a} - \sum_{k=0}^{\infty} \frac{1}{\left(k + b + \frac{1}{2}\right)^a} = \zeta(a, b) - \zeta\left(a, b + \frac{1}{2}\right).
\end{align*}
It is well known that \cite[(25.11.14)]{bib23}, \cite[(25.11.18)]{bib23}
\begin{equation}\label{fred}
\zeta(-n, x) = -\frac{1}{n+1}B_{n+1}(x),\quad \zp(0, x) = \ln\left(\frac{\Gamma(x)}{\sqrt{2\pi}}\right).
\end{equation}
By means of \eqref{fred}, we conclude the proofs of \eqref{bernz} and \eqref{zzerox}. The reflection formula \eqref{adamchik1} follows from \cite[(9)]{bib15}. It is well known that
\begin{equation}\label{rrfsa}
\sum_{k=0}^\infty \frac{(-1)^k}{k + z} = -\frac{1}{2}\left(\psi\left(\frac{z}{2}\right) - \psi\left(\frac{z+1}{2}\right)\right), \quad z \in \mathbb{C}\setminus -\mathbb{N}_0.
\end{equation}
Differentiating both sides of \eqref{rrfsa} with respect to $z$, $n-1$ times, we derive
\begin{equation}\label{rrfsb}
\sum_{k=0}^\infty \frac{(-1)^k}{(k + z)^n} = \frac{(-1)^{n}}{2^{n}(n-1)!}\left(\psi_{n-1}\left(\frac{z}{2}\right) - \psi_{n-1}\left(\frac{z+1}{2}\right)\right).
\end{equation}
We conclude the proof of \eqref{zpsi} from \eqref{rrfsb}. By virtue of two series \cite[(25.2.6), (25.2.11)]{bib23}
$$\sum_{k=1}^\infty \frac{(-1)^k \ln{k}}{k} = \gamma \ln{2} - \frac{1}{2}\ln^2{2}, \quad \sum_{k=1}^{\infty}\frac{(-1)^k\ln{k}}{k^n} = \frac{\ln{2}\zeta(n) + \left(2^{n-1} - 1\right)\zp(n)}{2^{n-1}}, \,\, n \geq 2,$$
\eqref{fonl} and \eqref{zpsi}, we conclude the proofs of  \eqref{prlf}.
\end{proof}
\end{lemma}
\section{The value of $\Delta_1$ and a generalization of $\Delta_n$}\label{secvardi}
In this section, we first illustrate the derivation of $\Delta_1$ from Glasser's result and then present an alternative proof for $\Delta_1$.
\subsection{Derivation of $\Delta_1$ from Glasser's result}
In this subsection, we show that $\Delta_1$ arises from Glasser's evaluation of $\beta^{\prime}(0)$. Glasser \cite[(A14)]{bib28} showed in 1972, that 
\begin{equation}\label{zref13}
\beta^{\prime}(0) = -2\ln\left(\frac{2\Gamma\left(\frac{3}{4}\right)}{\Gamma\left(\frac{1}{4}\right)}\right).
\end{equation}
Evaluating $\Delta_1$ by term-wise integration
\begin{equation}\label{zref12}
\Delta_1  = -\frac{\pi\gamma}{2} + 2\beta^{\prime}(1).
\end{equation}
The reflection formula for $\beta(s)$, for $s \in \mathbb{C}\setminus\{1\}$, is \cite[(A7)]{bib28}
\begin{equation}\label{zref11}
\beta(s) = \left(\frac{\pi}{2}\right)^{s-1}\cos\left(\frac{\pi s}{2}\right)\Gamma(1-s)\beta(1-s).
\end{equation}
Differentiating \eqref{zref11} at $s=0$, we derive
\begin{equation}\label{zref14}
\beta^{\prime}(0) = \frac{\gamma}{2} + \frac{1}{2}\ln\left(\frac{\pi}{2}\right) - \frac{2}{\pi}\beta^{\prime}(1).
\end{equation}
We infer from \eqref{zref13}, \eqref{zref12} and \eqref{zref14}
\begin{equation}\label{interm}
 -2\ln\left(\frac{2\Gamma\left(\frac{3}{4}\right)}{\Gamma\left(\frac{1}{4}\right)}\right) = \frac{\gamma}{2} + \frac{1}{2}\ln\left(\frac{\pi}{2}\right) - \frac{1}{\pi}\left(\Delta_1 + \frac{\pi\gamma}{2}\right).
\end{equation}
Expressing \eqref{interm} in terms of $\Delta_1$, we conclude the proof of \eqref{he2r}.
\subsection{A new method for evaluating $\Delta_1$}
We evaluate \eqref{he2r} with $\ln{x}$ replaced with $\ln{(x^2 + a^2)}$ for real values of $a$. The latter integral yields integrals \eqref{he2r} at $a=0$. This integral was proven differently via contour integration in Blagouchine's paper \cite[(4.1.1.1)]{bib1}. A generalization via contour integration can also be found in the same paper \cite[p.~80]{bib1}.
\begin{lemma}[Blagouchine 2014]\label{propf} Let $\Re\,b>0$. Then
\begin{align}
&\int_0^{\infty} \ln(x^2 + a^2) \sech(bx)\, \rmd{d}x  = \frac{2\pi}{b}\ln{\left(\frac{\sqrt{\frac{2\pi}{b}}\Gamma\left(\frac{b\lvert a \rvert }{2\pi} + \frac{3}{4}\right)}{\Gamma\left(\frac{b\lvert a \rvert }{2\pi} + \frac{1}{4}\right)}\right)}, \quad a \in \mathbb{R}. \label{propi}
\end{align}
\end{lemma}

\begin{proof}
We have
\begin{align}
&\int_0^{\infty} \frac{\ln(x^2 +  a^2)}{\cosh\left(\pi x\right)} \, \mathrm{d}x - \ln{a} = - \frac{2i}{\pi}\int_0^{\infty} \frac{\arctan\left(\mathrm{e}^{-\pi x}\right)}{\lvert a \rvert  - ix}\, \mathrm{d}x + \frac{2i}{\pi}\int_0^{\infty} \frac{\arctan\left(\mathrm{e}^{-\pi x}\right)}{\lvert a \rvert  + ix} \, \mathrm{d}x\nonumber
\\&= - \frac{2i}{\pi}\int_0^{\infty}\mathrm{e}^{-\lvert a \rvert t}\int_0^{\infty} \left(\mathrm{e}^{itx} - \mathrm{e}^{-itx}\right) \arctan\left(\mathrm{e}^{-\pi x}\right)\mathrm{d}x \, \mathrm{d}t\nonumber
\\&=  \frac{4}{\pi}\int_0^{\infty}\mathrm{e}^{-\lvert a \rvert t}\left( \frac{\pi}{4t} - \frac{1}{2t}\int_0^{\infty}\frac{\cos{\left(\frac{tx}{\pi}\right)}}{\cosh{x}}\, \mathrm{d}x\right) \mathrm{d}t = -\int_0^{1} \frac{z^{2\left|a\right|}}{\ln{z}}\left(\frac{\left(1 - z\right)^2}{z\left(1 + z^2\right)}\right)\, \mathrm{d}z\label{rzq2}
\\&= \int_0^{1}\int_0^1 \frac{z^{2\left|a\right| + p - 1}\left(1 - z\right)}{1 + z^2} \mathrm{d}z\, \mathrm{d}p = \frac{1}{2}\int_0^{1} \sum_{k=0}^{\infty} \left(-1\right)^k \left(\frac{1}{k + \frac{2\left|a\right| + p}{2}} - \frac{1}{k + \frac{2\left|a\right| + p + 1}{2}}\right)\mathrm{d}p \label{rzq13}
\\&= -\frac{1}{4}\int_0^{1} \left(\psi\left(\frac{2\left|a\right| + p}{4}\right) - \psi\left(\frac{2\left|a\right| + p+2}{4}\right) - \psi\left(\frac{2\left|a\right| + p + 1}{4}\right)  \right.\nonumber\\&\qquad\qquad\left. + \psi\left(\frac{2\left|a\right| + p + 3}{4}\right)\right)  \mathrm{d}p. \nonumber
\end{align}
Hence,
\begin{equation}
\int_0^{\infty} \frac{\ln(x^2 +  a^2)}{\cosh\left(\pi x\right)} \, \mathrm{d}x   = 2\ln\left(\frac{\sqrt{2}\Gamma\left(\frac{\lvert a \rvert }{2} + \frac{3}{4}\right)}{\Gamma\left(\frac{\lvert a \rvert }{2} + \frac{1}{4}\right)}\right).\label{rzq}
\end{equation}
We first show that $\int_0^\infty \cos{(\tau x)}\sech{x} \, \rmd{d}x = \pi/2 \sech\left(\pi \tau/2\right)$, for $\tau \in \mathbb{R}$, then we easily conclude the second part of \eqref{rzq2}. By writing the infinite sum on the second part of  \eqref{rzq13} in terms of the digamma function and integrating from $0$ to $1$, we finally have \eqref{rzq}. By making the substitution $x \mapsto bx/\pi$ in \eqref{rzq}, we conclude the proof of \eqref{propi}. To prove \eqref{he2r}, evaluate the limit of both sides of \eqref{rzq} as $a \to 0$.
\end{proof}
\subsection{A generalization of $\Delta_n$}
We generalize $\Delta_n$ as $\Delta_n(a, b)$, where $\Delta_n(1, 1)= \Delta_n$.
\begin{theorem}[Generalization of $\Delta_n$]\label{newthm} Let $n \in \mathbb{N}$, $a > 0$ and $\Re\, b > 0$. Then
\begin{align}
\Delta_n(a, b) &:= \int_0^\infty \ln{(ax)}\sech^n{(bx)} \, \rmd{d}x \nonumber
\\&\hspace{0.1cm}= \frac{2^{n-2}\Gamma^2\left(\frac{n}{2}\right)\ln\left(\frac{a}{b}\right)}{b(n-1)!} + \frac{2^{2n-1}}{b(n-1)!}  \sum_{k=2}^{n+1} \left(-\frac{1}{2}\right)^{k} P_{k-2}\left(n-1, \frac{n}{2}\right) \left[\zp\left(k-n, \frac{n}{4}\right)\right. \nonumber
\\&\hspace{0.6cm}- \zp\left(k-n, \frac{n+2}{4} \right) - \left(\gamma + \ln{4}\right)\left.\left(\zeta\left(k-n, \frac{n}{4}\right) - \zeta\left(k-n, \frac{n+2}{4}\right)\right)\right].\label{impm}
\end{align}
\end{theorem}
\begin{proof}
We first show that $\int_0^\infty \sech^n{x} \, \rmd{d}x = 2^{n-2}\Gamma^2\left(\frac{n}{2}\right)/(n-1)!$, then \eqref{impm} follows from Ripon's generalization of the $\ln$-$\ln$ form of $\Delta_n$ \cite[(2.1)]{bib16}.
\end{proof}
In the following lemma, we rediscover the limit $\lim_{a \to 1} \frac{\rmd{d}}{\rmd{d}a} 2^a\Phi\left(-1, a, \frac{1}{2}\right)$. A generalization of this appears in Blagouchine's paper \cite[p.~97, (b.2)]{bib1}.
\begin{lemma}\label{lem2}
We discover from Theorem \ref{newthm} at $a=b=1$, $n=1$, and \eqref{he2r}:
\begin{equation}
\lim_{n\to1}\left[\zp\left(n,  \frac{1}{4}\right) - \zp\left(n,   \frac{3}{4}\right)\right] = \lim_{a \to 1} \frac{\rmd{d}}{\rmd{d}a} 2^a\Phi\left(-1, a, \frac{1}{2}\right) = 2\pi\ln\left(\frac{\Gamma\left(\frac{3}{4}\right)\sqrt{2\pi}}{\Gamma\left(\frac{1}{4}\right)}\right) + \pi\left(\gamma + \ln{4}\right).\label{zpi}
\end{equation}

\end{lemma}
\begin{remark}
Implementing the $\Phi$-representation, \textsf{Maple~2022} evaluates \eqref{zpi} as $\simeq 5.12678$. Implementing the $\zeta$-representation, \textsf{Maple~2022} wrongly evaluates \eqref{zpi} as $\pi$. \textsf{Mathematica~13} evaluates correctly both the $\Phi$ and $\zeta$-representations.
\end{remark}
\begin{ex}\label{scex}
The following results are derived from Theorem \ref{newthm}:
$$\Delta_1(a, b) =  \frac{\pi}{b}\ln\left(\frac{\Gamma\left(\frac{3}{4}\right)\sqrt{2a\pi}}{\sqrt{b}\Gamma\left(\frac{1}{4}\right)}\right), \quad \Delta_2(a,b) =  -\frac{\gamma}{a} + \frac{1}{a}\ln\left(\frac{b\pi}{4a}\right),$$
$$\Delta_3(a,b)  = -\frac{2\textbf{G}}{b\pi} + \frac{\pi}{2b}\ln\left(\frac{\sqrt{2a\pi}\Gamma\left(\frac{3}{4}\right)}{\sqrt{b}\Gamma\left(\frac{1}{4}\right)}\right),\quad \Delta_4(a,b)= -\frac{2\gamma}{3b} - \frac{7\zeta(3)}{3b\pi^2} + \frac{2}{3b}\ln\left(\frac{a\pi}{4b}\right),$$
$$\Delta_5(a,b) =  -\frac{5\G}{3b\pi} + \frac{\pi}{24b} - \frac{\psi_3\left(\frac{1}{4}\right)}{192b\pi^3}+ \frac{3\pi}{8b}\ln\left(\frac{\sqrt{2a\pi}\Gamma\left(\frac{3}{4}\right)}{\sqrt{b}\Gamma\left(\frac{1}{4}\right)}\right),$$
$$\Delta_6(a,b) = -\frac{8\gamma}{15b} + \frac{8}{15b}\ln\left(\frac{a\pi}{4b}\right) - \frac{7\zeta(3)}{3b\pi^2} - \frac{31\zeta(5)}{5b\pi^4}.$$
\end{ex}
\begin{remark}
The relationship between $\Delta_1(a, b)$ and $\beta^{\prime}(1)$ is:
\begin{equation}
\Delta_1(a, b) = \frac{\pi}{2b}\left(\ln\left(\frac{a}{b}\right)-\gamma\right) + \frac{2}{b}\beta^{\prime}(1).
\end{equation}
\end{remark}
\section{A reduction formula and expressions for the integral sequences}\label{section4}
We show that relationship \eqref{red4} yields a reduction formula for $\Delta_n$. We also discover an expression for a new family of integrals by summing up $\chi_{k_1}$ from $k_1=l_1$ to $k_2$, $k_2=l_2$ to $k_3$ until $k_n=l_n$ to $n$, where $l_n \in \mathbb{N}$, for each $n \in \mathbb{N}$.

\subsection{A reduction formula for $\Delta_n$}
From a reduction formula for $\delta_n$, we derive a new reduction formula for $\Delta_n$.
\begin{proposition} \label{ssect3.7} Let $n\geq4$, where $n \in \mathbb{N}$. Then
\begin{align}
&\Delta_n = -\frac{4\textbf{\textup{G}}}{\pi(n-2)(n-1)} + \frac{n-2}{n-1}\Delta_{n-2}  + \frac{1}{(n-2)(n-1)}\sum_{k=1}^{n-3} \int_0^1\int_0^1\frac{\Gamma\left(\frac{a+b+k}{2}\right)}{k!} \nonumber \\ 
&\times\left[\sum_{j=0}^{k-1}\frac{2^{j-1}(k-j-1)!}{\Gamma\left(\frac{a+b+k-2j}{2}\right)}+\frac{1}{4\Gamma\left(\frac{a+b-k}{2}\right)}\left(\psi\left(\frac{a+b-k}{4}\right)- \psi\left(\frac{a+ b - k + 2}{4}\right)\right)\right]\rmd{d}a\, \rmd{d}b.\label{thefinal} 
\end{align}
\end{proposition}
\begin{proof} We derive a reduction formula for $\delta_n$ as:
\begin{align}
\delta_n &= -2^n\int_0^1\int_0^1 \frac{\Gamma\left(\frac{a+b+n}{2}\right)}{n!}\left[\sum_{k=0}^{n-1} \frac{2^{k-n-1}(n-k-1)!}{\Gamma\left(\frac{a+b+n-2k}{2}\right)}\right. \nonumber
\\&\qquad\qquad\quad \left.+\frac{1}{4\Gamma\left(\frac{a+b-n}{2}\right)}\left(\psi\left(\frac{a+b-n}{4}\right)- \psi\left(\frac{a+ b - n + 2}{4}\right)\right)\right]\rmd{d}a\, \rmd{d}b.\label{reftofinal}
\end{align}
By substituting $n=3$ in \eqref{eqn1.5}, we obtain the value of $\lambda_1$ using the values of $\Delta_1$ and $\Delta_3$ in Example \ref{scex}. Changing the dummy variable of the summation in \eqref{reftofinal} to $j$, substituting $n=k$ in \eqref{reftofinal}, and substituting the value of $\lambda_1$ in \eqref{red4}, we conclude the proof of \eqref{thefinal}.
\end{proof}

\subsection{Expressions for the integral sequences}
Taking account of \eqref{eqn1.5}, \eqref{eqn1.6}, \eqref{red2} and Theorem \ref{newthm}, we derive expressions for $\lambda_n$, $\delta_n$ and $\chi_n$. These expressions can be simplified with the following lemma.
\begin{lemma}\label{lemconj}Let $n \in \mathbb{N}$. Then
\begin{align}
&\frac{2^{2n+1}}{(n-1)!}\left[\frac{n^2}{16}\sum_{k=0}^{n-1} \left(-\frac{1}{2}\right)^k P_k\left(n-1, \frac{n}{2}\right) \left(\zeta\left(k-n+2, \frac{n}{4}\right) -\zeta\left(k-n+2, \frac{n+2}{4}\right)\right)\right.\nonumber
\\&- \left.\sum_{k=0}^{n+1} \left(-\frac{1}{2}\right)^k P_k\left(n+1, \frac{n+2}{2}\right) \left(\zeta\left(k-n, \frac{n+2}{4}\right) - \zeta\left(k-n, \frac{n+4}{4}\right)\right)\right] = 0.\label{conj1}
\end{align}
\end{lemma}
\begin{proof}
Integrating term-wise, one can easily show that
\begin{align}
\int_0^1 \frac{x^{n-1}}{(1+x^2)^n} \rmd{d}x = &\frac{2^{n-3}}{(n-1)!}\sum_{k=0}^{n-1} \left(-\frac{1}{2}\right)^k P_k\left(n-1, \frac{n}{2}\right) \nonumber
\\&\times\left(\zeta\left(k-n+2, \frac{n}{4}\right) -\zeta\left(k-n+2, \frac{n+2}{4}\right)\right). \label{conj1a}
\end{align}
We have by elementary integration
\begin{equation}\label{fanh10}
\int_0^1 \frac{x^{n-1}}{(1+x^2)^n} \rmd{d}x = \frac{\Gamma^2\left(\frac{n}{2}\right)}{4(n-1)!}.
\end{equation}
Equations \eqref{conj1a} and \eqref{fanh10} immediately yield the relationship
\begin{equation}\label{cojer}
\sum_{k=0}^{n-1} \left(-\frac{1}{2}\right)^k P_k\left(n-1, \frac{n}{2}\right) \left(\zeta\left(k-n+2, \frac{n}{4}\right) -\zeta\left(k-n+2, \frac{n+2}{4}\right)\right) = 2^{1-n}\Gamma^2\left(\frac{n}{2}\right).
\end{equation}
Substituting $n$ with $n + 2$ in \eqref{cojer}, and subtracting the resulting expression from the product of $n^2/16$ and \eqref{cojer}, \eqref{conj1} follows. 
\end{proof}
In the next theorem, we give new expressions for the integral sequences.
\begin{theorem}[New expressions for $\chi_n$, $\lambda_n$  and $\delta_n$]\label{conjecture}
The integral sequences $\chi_n$, $\lambda_n$  and $\delta_n$ can be expressed as follows:
\begin{align}
&\chi_{2n} = -\frac{4\G}{\pi} - \frac{2^{4n+1}}{(2n-1)!}  \sum_{k=1}^{n} \frac{1}{4^k} \left(n^2 P_{2k-2}\left(2n-1, n\right)  +  P_{2k}\left(2n+1, n+1\right)\right) \nonumber
\\&\qquad \times \left(\zp\left(-2n+2k, \frac{n}{2}\right) - \zp\left(-2n+2k, \frac{n+1}{2}\right)\right)  +\frac{2^{4n+1}n^{2n}\ln\frac{n}{2}}{4^n(2n-1)!}  \sum_{k=1}^{n}\frac{1}{n^{2k}} \nonumber
\\&\qquad \times P_{2k}\left(2n+1, n+1\right) - \frac{2^{4n+1}}{(2n-1)!}\left(\zp\left(-2n, \frac{n+1}{2}\right)  - \zp\left(-2n, \frac{n+2}{2}\right)\right),\label{conj2}
\end{align}
\begin{align}
&\chi_{2n+1} = -\frac{4\G}{\pi} + \frac{2^{4n-1}(2n + 1)^2}{(2n)!}\sum_{k=0}^n \frac{1}{4^k} P_{2k}\left(2n, \frac{2n+1}{2}\right) \nonumber
\\&\qquad\times\left(\zp\left(-2n+2k+1, \frac{2n+1}{4}\right) - \zp\left(-2n+2k+1, \frac{2n+3}{4}\right)\right) - \frac{2^{4n+3}}{(2n)!} \sum_{k=0}^{n+1} \frac{1}{4^k} \nonumber
\\&\qquad\times P_{2k}\left(2n+2, \frac{2n+3}{2}\right)\left(\zp\left(-2n+2k-1, \frac{2n+3}{4}\right) - \zp\left(-2n+2k-1, \frac{2n+5}{4}\right)\right), \label{conj3}
\end{align}

\begin{align}
&\lambda_{2n} = -\frac{2^{4n}}{n(2n-1)!}  \sum_{k=1}^{n} \frac{1}{4^k} \left(n^2 P_{2k-2}\left(2n-1, n\right)  + P_{2k}\left(2n+1, n+1\right)\right) \nonumber
\\&\qquad\times \left(\zp\left(-2n+2k, \frac{n}{2}\right) - \zp\left(-2n+2k, \frac{n+1}{2}\right)\right)  +\frac{2^{2n}n^{2n-1}\ln\frac{n}{2}}{(2n-1)!}  \sum_{k=1}^{n}\frac{1}{n^{2k}} \nonumber
\\&\qquad\times P_{2k}\left(2n+1, n+1\right) - \frac{2^{4n}}{n(2n-1)!}\left(\zp\left(-2n, \frac{n+1}{2}\right)  - \zp\left(-2n, \frac{n+2}{2}\right)\right), \label{conj4}
\end{align}
\begin{align}
&\lambda_{2n+1} = \frac{2^{4n-1}(2n + 1)}{(2n)!}\sum_{k=0}^n \frac{1}{4^k} P_{2k}\left(2n, \frac{2n+1}{2}\right)\nonumber
\\&\qquad\times\left(\zp\left(-2n+2k+1, \frac{2n+1}{4}\right) - \zp\left(-2n+2k+1, \frac{2n+3}{4}\right)\right) - \frac{2^{4n+3}}{(2n+1)!} \sum_{k=0}^{n+1} \frac{1}{4^k}  \nonumber
\\&\qquad\times P_{2k}\left(2n+2, \frac{2n+3}{2}\right)\left(\zp\left(-2n+2k-1, \frac{2n+3}{4}\right) - \zp\left(-2n+2k-1, \frac{2n+5}{4}\right)\right), \label{conj5}
\end{align}

\begin{align}
&\delta_n = \frac{2^{2n-1}(n+1)^2}{n!}  \sum_{k=0}^{n} \left(-\frac{1}{2}\right)^k P_k\left(n, \frac{n+1}{2}\right)\left(\zp\left(k-n+1, \frac{n+1}{4}\right)\right.\nonumber
\\&\qquad\left.- \zp\left(k-n+1, \frac{n+3}{4}\right)\right) -\frac{2^{2n+3}}{n!}  \sum_{k=0}^{n+2} \left(-\frac{1}{2}\right)^k P_k\left(n+2, \frac{n+3}{2}\right)  \nonumber
\\&\qquad\times \left(\zp\left(k-n-1, \frac{n+3}{4}\right) - \zp\left(k-n-1, \frac{n+5}{4}\right)\right) - \frac{2^{2n-3}n^2}{(n-1)!}  \sum_{k=0}^{n-1} \left(-\frac{1}{2}\right)^k\nonumber
\\&\qquad\times P_k\left(n-1, \frac{n}{2}\right)\left(\zp\left(k-n+2, \frac{n}{4}\right)- \zp\left(k-n+2, \frac{n+2}{4}\right)\right) + \frac{2^{2n+1}}{(n-1)!}  \nonumber
\\&\qquad\times  \sum_{k=0}^{n+1} \left(-\frac{1}{2}\right)^k P_k\left(n+1, \frac{n+2}{2}\right)   \left(\zp\left(k-n, \frac{n+2}{4}\right) - \zp\left(k-n, \frac{n+4}{4}\right)\right). \label{conj6}
\end{align}
\end{theorem}
\begin{proof}
Applying Lemma \ref{lemconj}, $\chi_n$ reduces to
\begin{align}
\chi_{n} = -\frac{4\G}{\pi} &+ \frac{2^{2n-3}n^2}{(n-1)!}  \sum_{k=0}^{n-1} \left(-\frac{1}{2}\right)^k P_k\left(n-1, \frac{n}{2}\right)\left(\zp\left(k-n+2, \frac{n}{4}\right)\right.\nonumber
\\&\left.- \zp\left(k-n+2, \frac{n+2}{4}\right)\right) -\frac{2^{2n+1}}{(n-1)!}  \sum_{k=0}^{n+1} \left(-\frac{1}{2}\right)^k P_k\left(n+1, \frac{n+2}{2}\right)  \nonumber
\\&\times \left(\zp\left(k-n, \frac{n+2}{4}\right) - \zp\left(k-n, \frac{n+4}{4}\right)\right). \label{conj2al}
\end{align}
Making the substitutions $n\mapsto 2n$, $n\mapsto 2n+1$ in \eqref{conj2al}, and taking account of Lemma \ref{lem0}, we conclude the proofs of \eqref{conj2} and \eqref{conj3}. Using the relationship between $\lambda_n$ and $\chi_n$ in \eqref{red2}, and using \eqref{conj2} and \eqref{conj3}, we conclude the proofs of \eqref{conj4} and \eqref{conj5}. Using the relationship between $\delta_n$ and $\lambda_n$ in \eqref{eqn1.6}, and using \eqref{conj4} and \eqref{conj5}, we conclude the proof of \eqref{conj6}.
\end{proof}

\begin{ex}\label{fex}The following results are derived from Theorem \ref{conjecture}. For $\chi_n$, one has:
$$\chi_1 = 0, \quad \chi_2 =  - \frac{4\textbf{G}}{\pi} + \frac{14\zeta(3)}{\pi^2}, \quad \chi_3 = -\frac{2\G}{\pi} - \frac{\pi}{2} + \frac{\psi_3\left(\frac{1}{4}\right)}{16\pi^3},$$
$$\chi_4 = - \frac{4\textbf{G}}{\pi} + \frac{28\zeta(3)}{3\pi^2} + \frac{124\zeta(5)}{\pi^4},\quad \chi_5 = -\frac{5\G}{2 \pi}-\frac{3 \pi}{4}+\frac{5 \psi_3\left(\frac{1}{4}\right)}{96 \pi^{3}}+\frac{\psi_5\left(\frac{1}{4}\right)}{768 \pi^{5}},$$
$$\chi_6 = - \frac{4\textbf{G}}{\pi} + \frac{112\zeta(3)}{15\pi^2} + \frac{124\zeta(5)}{\pi^4} + \frac{762\zeta(7)}{\pi^6}.$$ 
For $\delta_n$, one has:
$$\delta_1 =  -\frac{4\textbf{G}}{\pi} + \frac{14\zeta(3)}{\pi^2}, \quad \delta_2  = -\frac{14\zeta(3)}{\pi^2} + \frac{2\G}{\pi} - \frac{\pi}{2} + \frac{\psi_3\left(\frac{1}{4}\right)}{16\pi^3},$$
$$\delta_3 = -\frac{2\G}{\pi} + \frac{28\zeta(3)}{3\pi^2} + \frac{124\zeta(5)}{\pi^4} + \frac{\pi}{2} - \frac{\psi_3\left(\frac{1}{4}\right)}{16\pi^3},$$
$$\delta_4 = \frac{3\G}{2\pi} - \frac{28\zeta(3)}{3\pi^2} - \frac{124\zeta(5)}{\pi^4} - \frac{3\pi}{4} + \frac{5\psi_3\left(\frac{1}{4}\right)}{96\pi^3} + \frac{\psi_5\left(\frac{1}{4}\right)}{768\pi^5},$$
$$\delta_5 = -\frac{3\G}{2\pi} + \frac{112\zeta(3)}{15\pi^2} + \frac{124\zeta(5)}{\pi^4} + \frac{762\zeta(7)}{\pi^6} + \frac{3\pi}{4} - \frac{5\psi_3\left(\frac{1}{4}\right)}{96\pi^3} - \frac{\psi_5\left(\frac{1}{4}\right)}{768\pi^5},$$
$$\delta_6 = \frac{5\G}{4 \pi}  -\frac{15 \pi}{16}-  \frac{112\zeta(3)}{15\pi^2} - \frac{124\zeta(5)}{\pi^4} - \frac{762\zeta(7)}{\pi^6} +\frac{259\psi_3\left(\frac{1}{4}\right)}{5760 \pi^{3}}+\frac{7 \psi_5\left(\frac{1}{4}\right)}{4608 \pi^{5}} + \frac{\psi_7\left(\frac{1}{4}\right)}{92160 \pi^{7}}.$$\\
For $\lambda_n$, one has:
$$\lambda_1 = \frac{4\textbf{G}}{\pi}, \quad \lambda_2  = \frac{7\zeta(3)}{\pi^2}, \quad  \lambda_3 = \frac{2\G}{3\pi} - \frac{\pi}{6} + \frac{\psi_3\left(\frac{1}{4}\right)}{48\pi^3},$$
$$\lambda_4  = \frac{7\zeta(3)}{3\pi^2} + \frac{31\zeta(5)}{\pi^4}, \quad\lambda_5 = \frac{3\G}{10\pi} - \frac{3\pi}{20}+ \frac{\psi_3\left(\frac{1}{4}\right)}{96\pi^3} + \frac{\psi_5\left(\frac{1}{4}\right)}{3840\pi^5},$$
$$\lambda_6 = \frac{56\zeta(3)}{45\pi^2} + \frac{62\zeta(5)}{3\pi^4} + \frac{127\zeta(7)}{\pi^6}.$$
\end{ex}

\begin{remark}\label{rmk4.1}
By exploiting the fact that $\zeta(-2n)=0$, for all $n \in \mathbb{N}$, we readily derive the finite double sum from \eqref{cojer}:
\begin{equation}\label{thisbec}
\sum_{k=0}^{2n} \sum_{r=1}^{2n} (-1)^r r^{4n-2k} P_{2k}\left(4n+1, 2n+1\right) = \frac{((2n)!)^2}{2} - \frac{P_{4n}\left(4n+1, 2n+1\right)}{2}.
\end{equation}
We note that $s(n, 0) = 0$, $s(n,1)=(-1)^{n-1}(n-1)!$, and then obtain from the first part of \eqref{modff}
$$\frac{\Gamma(x)}{\Gamma(x-n-1)}\left(\psi(x) - \psi(x-n-1)\right) = \sum_{k=0}^{n} s(n+2, k+2) (k+1) x^k.$$
We derive from \eqref{lstr1}
$$P_n(n+1, x) = \sum_{j=0}^{n} s(n+2, j+2) (j +1) x^j.$$
Hence, we discover a functional equation for the signed generalized Stirling polynomials of the first kind, as:
\begin{equation}\label{functional}
P_n(n+1, x)  = \frac{\Gamma(x)}{\Gamma(x-n-1)}\left(\psi(x) - \psi(x-n-1)\right).
\end{equation}
By setting $n=2r+1$ in \eqref{yufo}, and $n=2r-1$, $x=(2r+1)/2$ in \eqref{functional}, the functional equation \eqref{functional} can be used to prove a particular case of Lemma \ref{lem0}. It can also be easily seen from \eqref{functional} that $P_{4n}\left(4n+1, 2n+1\right) = \left((2n)!\right)^2$. Hence, \eqref{thisbec} becomes
\begin{equation}
\sum_{k=0}^{2n} \sum_{r=1}^{2n} (-1)^r r^{4n-2k} P_{2k}\left(4n+1, 2n+1\right) = 0.
\end{equation}
We also derive from a generalized form of \eqref{conj1a}, a generalization involving the signed generalized Stirling polynomials of the first kind and the Hurwitz zeta function
\begin{align}
&\sum_{k=1}^{m+1} \left(-\frac{1}{2}\right)^{k-1} \left[P_{k-1} \left(m, \frac{p}{n}\right) \left(\zeta\left(k-m, \frac{p}{2n}\right) -\zeta\left(k-m, \frac{p+n}{2n}\right)\right) \right.\nonumber
\\&\quad\left.+ P_{k-1}\left(m, \frac{nm+n-p}{n}\right) \left(\zeta\left(k-m, \frac{nm+n-p}{2n}\right) -\zeta\left(k-m, \frac{nm+2n-p}{2n}\right)\right) \right]\nonumber
\\&\quad= 2^{1-m}\Gamma\left(\frac{p}{n}\right)\Gamma\left(\frac{nm+n-p}{n}\right),
\end{align}
where $m-\frac{p}{n}+1, \frac{p}{n}\in \mathbb{C}\setminus -\mathbb{N}_0$. \\

By summing up $\chi_{k_1}$ from $k_1=l_1$ to $k_2$, $k_2=l_2$ to $k_3$ until $k_n=l_n$ to $n$ for each $l_n \in \mathbb{N}$, we derive the family of integrals
\begin{align}
&\sum_{k_{n}=l_n}^{n} \cdots \sum_{k_2=l_2}^{k_3} \sum_{k_1=l_1}^{k_2} \chi_{k_1} = \sum_{k_{n}=l_n}^{n} \cdots \sum_{k_2=l_2}^{k_3} \sum_{k_1=l_1}^{k_2} \left[-\frac{4\G}{\pi} + \frac{2^{2k_1-3}k_1^2}{(k_1-1)!}  \sum_{m=0}^{k_1-1} \left(-\frac{1}{2}\right)^m\right.\nonumber
\\&\quad\times P_m\left(k_1-1, \frac{k_1}{2}\right)\left(\zp\left(m-k_1+2, \frac{k_1}{4}\right) - \zp\left(m-k_1+2, \frac{k_1+2}{4}\right)\right) -\frac{2^{2k_1+1}}{(k_1-1)!} \nonumber
\\&\quad\left.\times\sum_{m=0}^{k_1+1} \left(-\frac{1}{2}\right)^m P_m\left(k_1+1, \frac{k_1+2}{2}\right) \left(\zp\left(m-k_1, \frac{k_1+2}{4}\right) - \zp\left(m-k_1, \frac{k_1+4}{4}\right)\right)\right].\label{conj1al}
\end{align}
\end{remark}

\subsection{New integrals and generalizations from V\u{a}lean's integral}\label{secval}
V\u{a}lean \cite[\S 3.40]{bib9} presented two proofs for \eqref{2b} below. We shall not be consider those proofs. The integral in  \eqref{2b} shall be split into two convergent parts. Both parts are useful in evaluating $\Delta_2$, and it shall be seen that the two derived integrals give rise to a generalization of these integrals, which are not considered in \cite{bib9}.
\begin{proposition}[Integrals involving rational functions and the hyperbolic tangent function]\label{prop2}  Let $y \in \mathbb{R}$. Then
\begin{align} \label{2ai}& \int_0^\infty\left(\frac{x}{x^2 + y^2} - \frac{x}{x^2 + y^2}\tanh\left(\frac{\pi x}{2}\right)\right) \mathrm{d}x = \psi\left(\frac{1 + \lvert y \rvert }{2}\right) + \ln{\frac{2}{\lvert y \rvert }},\,\, y \neq 0, \\
\label{2aii}&\int_0^\infty\left(\frac{\tanh\left(\frac{\pi x}{2}\right)}{x} - \frac{x}{x^2 + y^2}\right)\mathrm{d}x = \gamma + \ln\left(2\lvert y \rvert \right), \quad y\neq 0,
\\ \label{2b}&\int_0^\infty \tanh\left( \frac{\pi x}{2}\right)\left(\frac{1}{x} - \frac{x}{x^2 + y^2}\right) \mathrm{d}x = -\psi\left(\frac{1}{2}\right) + \psi\left(\frac{1 + \lvert y \rvert }{2}\right).
\end{align}
\end{proposition}
\begin{proof} Recalling a product formula for the $\cosh$ function \cite[(4.36.2)]{bib23}
\begin{equation}\label{ldsksjdcss}
\cosh{x} := \prod_{k=0}^\infty \left(1 + \frac{4x^2}{\pi^2\left(2k+1\right)^2}\right),\quad x \in \mathbb{C},
\end{equation}
and evaluating the logarithmic derivative of both sides of \eqref{ldsksjdcss}, we obtain
\begin{equation}\label{ldsksjdcssdjdd}
\tanh{x} = 8x\sum_{k=0}^{\infty} \frac{1}{4x^2 + \pi^2\left(2k + 1\right)^2},\quad x \in \mathbb{C}.
\end{equation}
Introducing the variable $\epsilon_1 > 0$, replacing $x$ in \eqref{ldsksjdcssdjdd} with $\pi t/2$ and substituting in $\kappa_1(y)$, we derive
\begin{align}
\kappa_1(y) &= \lim_{N\to\infty}\frac{4}{\pi}\sum_{k = 0}^N \left(\frac{\lvert y \rvert }{\left(2k + 1\right)^2 - y^2}\arctan\left(\frac{\epsilon_1 N}{y}\right) - \frac{2k + 1}{\left(2k + 1\right)^2 - y^2}\arctan\left(\frac{\epsilon_1 N}{2k + 1}\right)\right)\nonumber
\\&\quad+ \frac{1}{2}\ln\left(\left(\frac{\epsilon_1 N}{y}\right)^2 + 1\right) = \psi\left(\frac{1 + \lvert y \rvert }{2}\right) - \ln{\lvert y \rvert } + \ln{\epsilon_1}.\label{eqn:3.19}
\end{align}
Evaluating the limit as $y\to \infty$ on both sides of \eqref{eqn:3.19} produces
\begin{align*}
\lim_{y \to \infty} \kappa_1(y)  &= \lim_{y \to \infty}\left(\psi\left(\frac{1 + \lvert y \rvert }{2}\right) - \ln{\lvert y \rvert } + \ln{\epsilon_1}\right),
\end{align*}
which implies that $\epsilon_1 = 2$. Therefore
\begin{align*}
\kappa_1(y)  &= \psi\left(\frac{1 + \lvert y \rvert }{2}\right) + \ln{\frac{2}{\lvert y \rvert }}. 
\end{align*}
Introducing the variable $\epsilon_2 > 0$, replacing $x$ in \eqref{ldsksjdcssdjdd} with $\pi t/2$ and substituting in $\kappa_2(y)$, we derive
\begin{align*}
\kappa_2(y)  &= \lim_{N \to \infty}\left(2H_{2N + 1} - 2\ln\left(2N + 1\right) - H_{N} + \ln{\frac{N}{\epsilon_2}}+ 2\ln\left(\frac{2N + 1}{N}\right)\right) +  \ln{\lvert y \rvert}
\\&= \gamma + \ln{\left(\frac{4\lvert y \rvert}{\epsilon_2}\right)}.
\end{align*} 
Evaluating the limit as $y \to 0$ on both sides of $\kappa(y)$ produces
\begin{align*}
0&= \lim_{y \to 0} \left(\gamma + \ln\left(\frac{4\lvert y \rvert }{\epsilon_2}\right) + \psi\left(\frac{1 + \lvert y \rvert }{2}\right) +  \ln{\frac{2}{\lvert y \rvert }}\right),
\end{align*}
which implies that $\epsilon_2 = 2$. Therefore $\kappa_2(y)  = \gamma + \ln\left(2\lvert y \rvert \right).$ Hence
\begin{align}\label{eqn:4}
\kappa(y) = -\psi\left(\frac{1}{2}\right) + \psi\left(\frac{1 + \lvert y \rvert }{2}\right).
\end{align}
\end{proof}
\begin{remark}
The value of $\Delta_2$ can be derived after showing from $\kappa_2(y)$ that
$$\kappa_2(y) =\lim_{N \to \infty}\left(\ln{N}\left(\tanh{N} - 1\right)- \frac{1}{2}\ln\left(1 + \frac{{\pi^2}{y^2}}{4N^2}\right)\right) + \ln\left(\frac{{\pi}{\lvert y \rvert }}{2}\right) - \Delta_2.$$
\end{remark}
\begin{proposition}[Generalization of Proposition \ref{prop2}]\label{prop3}
Let $y \in \mathbb{R}$ and $\Re\,a>0$. Then
\begin{align}
\label{31} \int_0^\infty \tanh\left(at\right)\left(\frac{1}{t} - \frac{t}{t^2 + y^2}\right) \mathrm{d}t &= -\psi\left(\frac{1}{2}\right) + \psi\left(\frac{1}{2} + \frac{a\lvert y \rvert }{\pi}\right),
\\
\label{32} \disp\int_0^\infty\frac{t}{t^2 + y^2}\left(1-\tanh\left(at\right)\right)\mathrm{d}t &= \psi\left(\frac{1}{2} + \frac{a\lvert y \rvert }{\pi}\right) + \ln\left(\frac{\pi}{a\lvert y \rvert }\right),\,\,y \neq 0,
\\
\label{34} \int_0^\infty\left(\frac{\tanh\left(at\right)}{t} - \frac{t}{t^2 + y^2}\right)\mathrm{d}t &= \gamma + \ln\left(\frac{4a\lvert y \rvert }{\pi}\right),\,\, y\neq 0.
\end{align}
\end{proposition}
\begin{proof}
To obtain the results in \eqref{31}, \eqref{32}, and \eqref{34}; replace $\lvert y \rvert $ with $2a\lvert y \rvert /\pi$ in \eqref{2ai}, \eqref{2aii}  and \eqref{2b}, and make the substitution $x=2at/\pi$.
\end{proof}

\begin{proposition}[Integrals involving the arctangent and hyperbolic tangent functions]\label{prop4}
Let $y \in \mathbb{R}\setminus\{0\}$ and $\Re\,a > 0$. Then\normalfont
\begin{align}
&\int_0^\infty\left(1 - \tanh\left(a t\right)\right)\arctan{\frac{\lvert y \rvert }{t}}\,\mathrm{d}t =\frac{\pi}{a}\ln{\left(\frac{\Gamma\left(\frac{1}{2} + \frac{a\lvert y \rvert }{\pi}\right)}{\sqrt{\pi}}\right)} + \lvert y \rvert \left(\ln{\left(\frac{\pi}{a\lvert y \rvert }\right)} + 1\right),\label{youpart}\\
&\int_0^\infty \frac{\ln{\left(\mathrm{e}^{at}\sech(at)\right)}}{t^2 + y^2} \, \rmd{d}t = \frac{\pi}{\lvert y \rvert }\ln{\left(\frac{\Gamma\left(\frac{1}{2} + \frac{a\lvert y \rvert }{\pi}\right)}{\sqrt{\pi}}\right)}+ a\left(\ln\left(\frac{\pi}{a\lvert y \rvert }\right) + 1\right).\label{younko}
\end{align}
\end{proposition}

\begin{proposition}[New integral representations for $\gamma$]\label{prop4w}
Let $y \in \mathbb{R}\setminus\{0\}$ and $\Re\,a > 0$. Then\normalfont
\begin{align}
\gamma &= -\ln{4} - \frac{\pi}{a\lvert y \rvert }\ln{\left(\frac{\Gamma\left(\frac{1}{2} + \frac{a\lvert y \rvert }{\pi}\right)}{\sqrt{\pi}}\right)} + \int_0^\infty \tanh{(at)}\left(\frac{1}{t} - \frac{1}{\lvert y \rvert }\arctan{\frac{\lvert y \rvert }{t}}\right) \rmd{d}t \label{excepb}
\\&= -\ln{4} - \frac{\pi}{a\lvert y \rvert }\ln{\left(\frac{\Gamma\left(\frac{1}{2} + \frac{a\lvert y \rvert }{\pi}\right)}{\sqrt{\pi}}\right)} + \frac{y^2}{a}\int_0^\infty \frac{\ln{\cosh\left(a t\right)}}{t^2\left(t^2 + y^2\right)} \mathrm{d}t
\\&= 1 - \ln\left(\frac{4a\lvert y \rvert }{\pi}\right) + \int_0^\infty\left(\frac{\ln{\cosh\left(a t\right)}}{at^2} - \frac{t}{t^2 + y^2}\right)\mathrm{d}t
\\&= 2-\ln\left(\frac{4a\lvert y \rvert }{\pi}\right) + \int_0^\infty \left(\frac{\ln{\cosh\left(a t\right)}}{at^2} - \frac{1}{\lvert y \rvert }\arctan{\frac{\lvert y \rvert }{t}}\right)\mathrm{d}t \label{lastf}
\\&= 1 - \ln\left(\frac{4a\lvert y \rvert }{\pi}\right) +\int_0^\infty \left(\frac{\tanh\left(a t\right)}{t} - \frac{1}{\lvert y \rvert }\arctan{\frac{\lvert y \rvert }{t}}\right)\mathrm{d}t\label{euu}.
\end{align}
\end{proposition}
\begin{proof}
To prove the results in Propositions \ref{prop4} and \ref{prop4w}; integrate both sides of each of the results in Proposition \ref{prop3} over $\left(0, a\right)$ with respect to $a$ and over $\left(0, \lvert y \rvert \right)$ with respect to $\lvert y \rvert $.
\end{proof}
\begin{ex} \label{thex} 
We have the following integral representations for $\gamma$.
\begin{align}
\gamma &= -\ln{2} + \int_0^\infty \tanh{t}\left(\frac{1}{t} - \frac{1}{\pi}\arctan{\frac{\pi}{t}}\right) \rmd{d}t = -\ln{2} + \int_0^\infty \frac{\pi^2\ln{\cosh{t}}}{t^4 + \pi^2t^2} \mathrm{d}t \nonumber
\\&=  2 +\int_0^\infty \left(\frac{\ln{\cosh\left(\pi t\right)}}{\pi t^2} - 4\arctan{\frac{1}{4t}}\right)\mathrm{d}t =  1 +\int_0^\infty \left(\frac{\tanh\left(\pi t\right)}{t} - 4\arctan{\frac{1}{4t}}\right)\mathrm{d}t\nonumber
\\\label{7we}&=  1 + \int_0^\infty\left(\frac{\ln{\cosh\left(\pi t\right)}}{\pi t^2} - \frac{16t}{16t^2 + 1}\right)\mathrm{d}t,
\end{align}
\end{ex}

\begin{ex} \label{thex2} 
We have the following integral representations for $\ln{2}$ and $\pi$.
$$\ln{2} = 1 - \int_0^\infty \frac{\ln\left(\mathrm{e}^t\sech{t}\right)}{t^2 + \pi^2} \, \rmd{d}t,$$
$$\pi := \frac{1}{1 - \ln{2}}\int_0^\infty\left(1 - \tanh{t}\right)\arctan{\frac{\pi}{t}}\,\mathrm{d}t.$$
\end{ex}
\begin{remark}
The integral representations are symmetric for $y$ and $a > 0$. It is obvious that \eqref{7we} are simpler expressions of $\gamma$ in Propositions \ref{prop4} and \ref{prop4w}. However, none of \eqref{youpart}, \eqref{excepb}, \eqref{lastf}, \eqref{euu} appear in Gradshteyn and Ryzhik \cite[(4.561), p.~605]{bib25}, and none of \eqref{excepb}--\eqref{euu} appear in Gradshteyn and Ryzhik \cite[(8.367), p.~905]{bib25}.
\end{remark}
\begin{remark}\label{rmk:lamb}
It can easily be deduced from \eqref{32} that
\begin{equation}
\int_0^\infty \ln(t^2 + y^2) \sech^2(at)\, \rmd{d}t = \frac{2}{a}\left(\psi\left(\frac{1}{2} + \frac{a\lvert y \rvert }{\pi}\right) + \ln{\frac{\pi}{a}}\right),
\end{equation}
where $\Re\,a > 0$ and $y \in \mathbb{R}$.
\end{remark}

\section*{Acknowledgment}
I acknowledge and appreciate the efforts of the Spirit of Ramanujan (SOR) STEM Talent Initiative, for providing me with the computational tools that I used in verifying my results. I thank Murphy E. Egwe, Howard S. Cohl, M. Lawrence Glasser, and Ken Ono for valuable discussions.
\section*{Funding}
The author did not receive funding from any organization for the submitted work.
\bibliographystyle{unsrt}
\bibliography{Abdulhafeez_A_Abdulsalam}
\end{document}